\input amstex 
\documentstyle{amsppt}
\input bull-ppt
\keyedby{bull491/paz}

\topmatter
\cvol{31}
\cvolyear{1994}
\cmonth{July}
\cyear{1994}
\cvolno{1}
\cpgs{39-43}
\ratitle
\title Closed ideals of the algebra of\\absolutely 
convergent
Taylor series\endtitle
\author J. Esterle, E. Strouse and F. Zouakia\endauthor
\shortauthor{J. Esterle, E. Strouse, and F. Zouakia}
\shorttitle{Algebra of absolutely convergent Taylor series}
\address ({\rm J. Esterle and E. Strouse})
UFR de Mathematiques et Informatique,
Universit\'e de Bordeaux I, 351 cours de la Liberation,
33405 Talence, France\endaddress
\ml\nofrills{\it E-mail address}{\rm, E. Strouse:}\enspace
\tt strouse\@ecole.ceremab.u-bordeaux.fr\endml
\address ({\rm F. Zouakia})
Department de Mathematiques, ENS Takkadoum, BP5118 Rabat, 
Morocco\endaddress
\date June 25, 1993\enddate
\abstract Let $\Gamma$ be the unit circle,
$A(\Gamma)$ the Wiener algebra of continuous functions
whose series of Fourier coefficients are absolutely 
convergent,
and $A^+$ the subalgebra of $A(\Gamma)$ of functions whose 
negative
coefficients are zero. If $I$ is a closed ideal of $A^+$, 
we denote
by $S_I$ the greatest common divisor of the inner factors 
of the
nonzero elements of $I$ and by $I^A$ the closed ideal 
generated by
$I$ in $A(\Gamma)$. It was conjectured that the equality 
$I^A=
S_I H^\infty \cap I^A$ holds for every closed ideal $I$. 
We exhibit
a large class $\scr F$ of perfect subsets of $\Gamma$, 
including
the triadic Cantor set, such that the above equality holds 
whenever
$h(I)\cap\Gamma\in\scr F$. We also give counterexamples to 
the
conjecture.\endabstract
\endtopmatter

\document
\heading 1. Introduction\endheading
Let $D$ be the open unit disk, let $\Gamma$ be the unit 
circle, and let
$H^\infty=H^\infty(D)$ be the algebra of bounded 
holomorphic functions on
$D$. Let $A^+$ be the algebra of absolutely convergent 
Taylor series,
i.e., the algebra of analytic functions $f\: D\to \Bbb C$ 
such that
$\|f\|_1=\sum^\infty_{n=0} \tfrac{|f^{(n)}(0)|}{n!} <+
\infty$. Clearly, $A^+$
is a subalgebra of the disc algebra $\scr U(D)$ consisting 
of functions
continuous on $\overline{D}$ and holomorphic on $D$. Also, 
if we
denote by $A(\Gamma)$ the algebra of absolutely convergent 
Fourier series,
equipped with the norm $\|f\|_1=\sum_{n\in{\Bbb Z}}|\hat 
f(n)|$, and we identify
$f\in A^+$ with $f_{|\Gamma}$, we can write $A^+=\{f\in 
A(\Gamma): \hat f(n)=0
(n<0)\}$.
\par
Let $I\ne\{0\}$ be a closed ideal of $A^+$, let $S_I$ be 
the greatest common
divisor of the inner factors of all nonzero elements of 
$I, I^A$ the closed
ideal of $A(\Gamma)$ generated by $I$, and 
$h(I)=\{z\in\overline{D}
\: f(z)=0 (f\in I)\}$. If $E\subseteq\Gamma$ is closed, 
set $I^+(E)=\{f
\in A^+: f_{|E} \equiv 0\}$. It was proved long ago by 
Carleson \cite{3} that
$I^+(E)=\{0\}$ for certain closed sets $E$ of measure 
zero. The structure of
those closed ideals $I$ of $A^+$ such that $h(I)$ is 
finite or countable was
described in 1972 by Kahane \cite{12} and Bennett-Gilbert 
\cite2. In this case,
$I=I^+(h(I)\cap\Gamma)\cap S_I\bfcdot H^\infty$. 
Bennett-Gilbert conjectured in
\cite 2 that, in general,
\thm{Conjecture 1}
$I=I^A\cap S_I\bfcdot H^\infty$.
\ethm
\par
(This conjecture was also quoted by Kahane \cite{12}.) 
Similar conjectures have
been verified by Beurling-Rudin \cite{19} for $\scr U(D)$, 
by Taylor-Williams
\cite{20} for $A^\infty(D)=\{f\in\scr U(D): 
f^{(n)}\in\scr{U}(D)
\ (n\ge1)\}$, by Korenblum \cite{16} for 
$A^p(D)=\{f\in\scr U(D):
f^{(n)}\in\scr U(D)\ (n\le p)\}$, and by Matheson 
\cite{17} for
$$
\Lambda_\alpha(D)=\{f\in\scr U(D):\lim
\frac{|f(e^{it_1})-f(e^{it_2})|}{|t_1-t_2|^\alpha}=0
\text{ uniformly as }|t_1-t_2|\to0)\}.
$$

The purpose of this note is to describe recent progress 
concerning closed
ideals of $A^+$. 

The main results can be described as follows:
\thm{Theorem 1.1 \rm\cite6}
Let $p$ be an integer $\ge3$.
If $h(I)\cap\Gamma\subseteq E_{1/p}$, the perfect 
symmetric set of constant
ratio $\tfrac1p$, then $I$ satisfies Conjecture \rm1.
\ethm
\thm{Theorem 1.2 \rm\cite7}
There exists a Kronecker set $E$ and a closed ideal $I$ of 
$A^+$ such that
$S_I=1$, $h(I)=E$, which does not satisfy Conjecture \rm1.
\ethm
\par
Results for $L^1(\Bbb R^+)$ which concern ideals whose 
hull is countable and
are analogous to the results of Kahane and Bennett-Gilbert 
have been obtained
by Nyman \cite{18} and Gurarii \cite{9, 10}. By using 
transfer methods due to
Hedenmalm \cite{11}, El Fallah \cite4 has proved certain 
versions of Theorem
1.1 and Theorem 1.2 for $L^1(\Bbb R^+)$.
\heading 2. Division ideals\endheading
If $I\ne\{0\}$ is a closed ideal of $A^+$ and if $f\in A^+
$, we denote by
$I(f)=\{g\in A^+: fg\in I\}$ the {\it division ideal\/} 
associated with $f$.
We will say that a closed subset $E\subseteq\Gamma$ is a 
Carleson set if
$$
\int^\pi_{-\pi} 
\log\left(\frac{2}{\roman{dist}(e^{it},E)}\right)\,dt<
\infty.
$$
\par
This condition is necessary and sufficient for the 
existence of a nonzero $f\in
\Lambda_\alpha$ vanishing on $E$. It is also necessary and 
sufficient for the
existence of an outer $f\in A^\infty(D)$ such that $f$ 
vanishes exactly on
$E$ and $f^{(n)}_{|E}=0$ $(n\ge1)$; see \cite{3, 20}.
\par
Some improvements of the methods used to discuss closed 
ideals of $A^\infty
(D)$ and $A^+(D)$ lead to the following result (one must 
circumvent
the fact that the outer part of an element of $A^+$ does 
not necessarily belong
to $A^+)$:
\thm{Theorem 2.1}
{\rm(i)} If $f\in A^+\cap S_I H^\infty$, then 
$h(I(f))\subseteq
h(I)\cap\Gamma$. Also, the positive singular measure which 
defines the inner
factor of $I(f)$ is nonatomic and vanishes on all Carleson 
sets.
\par
{\rm(ii)} If, further, $f$ vanishes on $h(I)\cap\Gamma$, 
then $h(I(f))$ is a
perfect subset of $\Gamma$.
\ethm
\par
We note that Theorem 2.1(ii) contains the results of 
Kahane and
Bennett-Gilbert, for, if $h(I)$ is countable, it implies 
that $h(I(f))=
\varnothing$.
\par
The results and methods of \cite{1, 17, 20} lead, also, to 
the following
information:
\thm{Theorem 2.2}
Let $E$ be a Carleson set, and denote by $J^+_0(E)$ the 
closure in $A^+$ of the
set of elements of $A^\infty(D)$ vanishing on $E$ with all 
their
derivatives. Then\,\RM:
\roster
\item"(i)" \<$J^+_0(E)\subseteq I$ for every closed ideal 
$I$ of $A^+$
such that $S_I=1$ and $h(I)\subseteq E$.
\item"(ii)" If $\alpha>\tfrac12$ and 
$f\in\Lambda_\alpha(D)\cap I^+
(E)$, then $f\in J^+_0(E)$.
\endroster
\ethm
\heading 3. When the conjecture works\endheading
The following lemma, related to the Katznelson-Tzafriri 
theorem for
contractions \cite{5, 14}, is the key to our positive 
results concerning the
Bennett-Gilbert conjecture. The $w^\star$ topology 
discussed below is defined
by considering $A^+$ as the dual of $c_0$.
\thm{Lemma 3.1}
Let $I$ be a closed ideal of $A^+$, and let $f\in I^A\cap 
A^+$. Then $I(f)$
is $w^\star$-closed.
\ethm
\par
Using the results of \S2, we obtain:
\thm{Theorem 3.2 \rm\cite6}
Let $E\subseteq\Gamma$ be a Carleson set. If $J^+_0(E)$ is 
$w^\star$-dense in
$A^+$, then $I=I^A\cap S_I H^\infty$ for every closed 
ideal of $A^+$ such that
the perfect part of $h(I)\cap\Gamma$ is contained in $E$.
\ethm
\par
Thus, Theorem 1.1 follows from the fact that $J^+_0 
(E_{1/p})$ is indeed
$w^\star$-dense in $A^+$.
\heading 4. When the conjecture fails\endheading
If $K$ is a Helson set, then $I^+(K)$ is $w^\star$-dense 
in $A^+$ (this
observation is an extension of \cite{8, Theorem 4.5.2}). 
Also, if $K$ is a
Kronecker set and a Carleson set, then $J^+_0(K)$ has no 
inner factor; so the
closed ideal generated by $J^+_0 (K)$ consists of all 
functions of $A(\Gamma)$
vanishing on $K$, since Kronecker sets satisfy synthesis 
\cite{21}. So, if the
ideal $J^+_0(K)$ satisfies the Bennett-Gilbert conjecture, 
we must have $J^+_0
(K)=I^+(K)$. A construction which has some relation with 
Kaufman's construction
of a Helson set of multiplicity \cite{15} gives the 
following result:
\thm{Theorem 4.1 \rm\cite7}
Let $E$ be a set of multiplicity. Then there exists a 
nonzero distribution
$\mu$ whose support is a Kronecker subset of $E$ such that 
$\hat\mu(n)\to0$ as
$n\to-\infty$.
\ethm
\par
Theorem 1.2 follows from Theorem 4.1 applied to a Carleson 
set of multiplicity
(for example, the perfect set $E_\xi$ when $\tfrac1\xi$ is 
not a Pisot number).
In this case, if $K$ is the support of the distribution 
given by Theorem 4.1,
we have $J^+_0(K)$ properly contained in $I^+(K)$, and it 
is even possible to
show that there are $2^{\aleph_0}$ distinct closed ideals 
between $J^+_0(K)$ and
$I^+(K)$, ideals which would be equal if the 
Bennett-Gilbert conjecture were
true.
\heading 5. Applications\endheading
The positive results about the Bennett-Gilbert conjecture 
give ``strong
uniqueness properties'' of some closed subsets of 
$\Gamma$. For example, it
follows from Theorem 2.1 that any distribution $S$ 
supported by $E_{1/p}$
such that $\widehat S(n)\to 0$ as $n\to\infty$ must be the 
zero distribution.
Some stronger results involving hyperdistributions can be 
found in \cite7. We
present here an application of Theorem 1.1 to operator 
theory (an extension of
the Beurling-Pollard method \cite{13, p.\ 61} is involved 
in the proof).
\thm{Theorem 5.1}
Let $T$ be a contraction on a Banach space. If 
$\roman{Sp}\,T\subseteq E_{1/p}$
and if
$$
\limsup_{n\to\infty}
\frac{\log^+\|T^{-n}\|}{n^\alpha}<+\infty\quad\text{where }
\alpha<\frac{\log p-\log2}{2\log p-\log 2},
$$
then $\sup_{n\ge1} \|T^{-n}\|<+\infty$.
\ethm
\par
If we add to the hypotheses of Theorem 5.1 the assumption 
that $\roman{Sp}\,T$
is a Dirichlet set, we can conclude that $T$ is an 
isometry. 
Analogous results hold
 for all $E_\xi$ with $\xi\in(0,\tfrac12)$ if we consider 
only
Hilbert spaces. On the other hand, if $E$ is a set of 
multiplicity (for
example, $E=E_\xi$ when $\tfrac1\xi$ is not a Pisot 
number), then there exist
contractions $T$ such that $\roman{Sp}\, T\subseteq E$ and 
$\|T^{-n}\|$ goes to
infinity arbitrarily slowly as $n$ goes to infinity.
\par
Zarrabi \cite{22} proved that if $E\subseteq\Gamma$ is 
countable, every
contraction on a Banach space $X$ such that 
$\roman{Sp}\,T\subseteq X$ and
$\log^+
\|T^{-n}\|/n^{1/2} \to0$ as $n\to\infty$ is an isometry; 
but this property
holds exclusively for countable sets, even if we suppose 
that $X$ is a Hilbert
space.
\par
It follows from unpublished computations
by M. Zarrabi, M. Rajoelina, and the first
author that the constant $\frac{\log p-\log 2}{2\log 
p-\log2}$
in Theorem 5.1 is the best possible.

\enddocument